\documentclass[11pt]{article}

\usepackage{amsmath,amsfonts,amssymb,amsthm,bbm,latexsym,mathrsfs}
\usepackage{graphicx,color,epsfig,fancyhdr,dsfont, ulem}
\usepackage{enumerate}
\usepackage{hyperref}
\usepackage{indentfirst}
\usepackage{amsmath,amscd}
\usepackage{caption}
\usepackage{graphicx, subfig}
\usepackage[all]{xy}
\usepackage[dvipsnames]{xcolor}
\title{On minimal flows of commutative $p$-adic groups}
\date{\today}

\author{Ningyuan Yao \\
  \text{Fudan University} \and Zhentao Zhang\\
  \text{Fudan University}}
\date{\today}

\newtheorem{theorem}{Theorem}[section]

\newtheorem{definition}[theorem]{Definition}
\newtheorem{rmk}[theorem]{Remark}
\newtheorem{lemma}[theorem]{Lemma}
\newtheorem{coro}[theorem]{Corollary}
\newtheorem{fact}[theorem]{Fact}

\newtheorem*{ques}{Newelski's Question}
\newtheorem{Conj}{Conjecture}
\newtheorem{notation}[theorem]{Notation}

\newtheorem*{acknowledgment}{Acknowledgments}

\newcommand{\tp}{\mathrm{tp}}
\newcommand{\Th}{\mathrm{Th}}
\newcommand{\cl}{\mathrm{cl}}
\newcommand{\Gen}{\mathrm{Gen}}
\newcommand{\ext}{\mathrm{ext}}
\newcommand{\dcl}{\mathrm{dcl}}
\newcommand{\Aut}{\mathrm{Aut}}
\newcommand{\Q}{\mathbb{Q}_p}
\newcommand{\Z}{\mathbb{Z}}
\newcommand{\M}{\mathbb{M}}

\newcommand{\N}{\mathbb{N}}

\newcommand{\Ga}{\mathbb{G}_\text{a}}
\newcommand{\Gm}{\mathbb{G}_\text{m}}

\newcommand{\pCF}{p\mathrm{CF}}
\newcommand{\WGen}{\mathrm{WGen}}

\newcommand{\WG}{\mathrm{WG}}
\newcommand{\AP}{\mathrm{AP}}

\newcommand{\SO}{\mathrm{SO}}

\newcommand{\st}{\mathrm{st}}
\newcommand{\I}{\mathcal{I}}
\newcommand{\J}{\mathcal{J}}

\newcommand{\dfg}{\mathrm{dfg}}
\newcommand{\fsg}{\mathrm{fsg}}
\newcommand{\NIP}{\mathrm{NIP}}
\newcommand{\img}{\mathrm{Im}}
\newcommand{\id}{\mathrm{id}}

\newcommand{\sq}{\subseteq}

\newcommand{\rd}{\upharpoonright}

\begin{document}

\maketitle

\begin{abstract}
We study the definable topological dynamics $(G,S_G(M))$ of a definable group acting on its type space, where $M$ is a structure and $G$ is a group definable in $M$. In \cite{Newelski-I}, Newelski raised a question of whether weakly generic types coincide with almost periodic types in definable topological dynamics. In \cite{YZ-Sta}, we introduced the notion of stationarity (see  Definition \ref{def-of-stationary}), showing the answer is positive when $G$ is a stationary definably amenable group definable over the field of $p$-adic numbers or an $o$-minimal expansion of real closed field.

In this paper, we continue with the work of \cite{YZ-Sta},  focusing on the case where $G$ is a commutative groups definable over the field of $p$-adic numbers,  and showing that weakly generic types coincide with  almost periodic types if and only if either $G$ has definable $f$-generics or $G$ is stationary.

\end{abstract}

\section{Introduction}

In \cite{Newelski-I}, Newelski introduced a connection between topological dynamics and model theory and now called definable topological
dynamics. Definable topological dynamics studies the action of $G(M)$, a  group defined in some model $M$, on $S_G(M)$, the space of types concentrating on $G$ over $M$ and tries to connect the concepts of such topological dynamics to the model theoretical ones. The problem raised in \cite{Newelski-I} of whether almost periodic types coincide with weakly generic types is in this topic.

The notion of ``weakly generic'' introduced by Newelski in \cite{Newelski-I} is a suitable generalization of ``generic'' from stable environments to unstable ones. We say that a definable set $X\sq G(M)$ is \emph{weakly generic} if there is a non-generic definable set $Y$ such that $X\cup Y$ is generic, where a definable set is \emph{generic} if its finitely many  translates cover the whole group. We say that a type $p\in S_G(M)$ is weakly generic if every definable set in it is weakly generic. We denote the space of weakly generic types in $S_G(M)$ by $\WGen(S_G(M))$.

The notion of ``almost periodic'' comes from topological dynamics. We say that a type $p\in S_G(M)$ is \emph{almost periodic} if the closure of its $G(M)$-orbit is a minimal subflow of $S_G(M)$. By $\AP(S_G(M))$, we denote the space of almost periodic types in $S_G(M)$.

Newelski proved in \cite{Newelski-I} that $\AP(S_G(M))$ coincides with the closure of $\WGen(S_G(M))$, and if the generic types exist, then $\AP(S_G(M))=\WGen(S_G(M))$. An example which has IP while is not simple was given in \cite{Newelski-I}  to show that the two classes differ. Newelski asked in \cite{Newelski-I} that
\begin{ques}\label{ques}
Assume that $G$ is definable in an $o$-minimal or even just an $\NIP$ structure. Is every (global) weakly generic type on $G$ almost periodic?
\end{ques}
In \cite{CS-amenable-$NIP$-group}, the Newelski's Question is restated in the special case of definably amenable groups in $\NIP$ theories. When $M$ is an $o$-minimal expansion of a real closed field and $G$ is a definably amenable group definable over $M$, Pillay and Yao proved in \cite{Pillay-Yao-mini-flow} that weakly generics coincide with almost periodics when the torsion free part of $G$ has dimension one. They also gave a counter-example when $(\mathbb R,+)^2\times \SO(2,\mathbb R)$ to show that the set of almost periodics is a proper subset of weakly generics.

Recall that a group $G$ has finitely satisfiable generics (fsg) or definable
$f$-generics (dfg) if there is a global type $p$ on $G$ and a small model $M$ such
that every left translate of $p$ is finitely satisfiable in $M$ or definable over $M$,
respectively.

In \cite{YZ-Sta}, we considered the Newelski's Question where $G$ is a  group definable over an $o$-minimal structure or the field $\Q$ of $p$-adic numbers and admitting a ``dfg-fsg decomposition'', namely, admitting a $\Q$-definable short exact sequence
\begin{equation}\label{dfg-fsg-decom}
    1\rightarrow H\rightarrow G\rightarrow_\pi C\rightarrow 1
\end{equation}
where $H$ is a  $\dfg$  group 
and $C$ is a  $\fsg$ group.
Note that a group definable in an $o$-minimal structure is definably amenable iff it has such decomposition (see \cite{CP-connected-component}). Let $M$ be either an $o$-minimal structure or the field of $p$-adics, we call a group definable over $M$ is stationary if every weakly generic type over $M$ has a unique global weakly generic extension. Let $G$ be a group definable over $M$ and satisfy the ``dfg-fsg decomposition'' as given in (\ref{dfg-fsg-decom}), we showed in \cite{YZ-Sta} that the almost periodics coincide  with the weakly generics if either  $G$ has $\dfg$  or $G$ is stationary, which extends the result of \cite{Pillay-Yao-mini-flow}. Also, we proved that $G$ is stationary if and only if $H$ is stationary.

In this article, we focus on Newelski's Question where $G$ is a commutative  group definable over $\Q$.  The advantage of working in the commutative case is that such groups admit the ``$\dfg$-$\fsg$ decomposition'' as given in (\ref{dfg-fsg-decom}) (see \cite{JY-abelian-pcf-group}). Hence, some arguments in \cite{YZ-Sta} works in the paper. Also, we can always assume $\dim(C)\geq 1$, because Yao proved in \cite{Yao-tri} that Newelski's question has a positive answer when $G=H$.

Our first result give a description for stationary commutative $\dfg$ groups over $\Q$:
\begin{theorem}\label{sta-no-Ga2}
Let $H$ be a commutative $\dfg$ group definable over $\Q$.  Then
\begin{enumerate}[(i)]
\item 
$H$ has a finite index subgroup $E$ and a finite subgroup $F$ such that $E/F$ is definably isomorphic to a finite index subgroup of $\Ga^s\times\Gm^t$ for some $s,t\in\mathbb{N}$, where $\Ga$ and $\Gm$ are the additive group and multiplicative group, respectively.
    \item Let $H$ be as in part (i),
    then $H$ is stationary iff   $s\leq 1$.
\end{enumerate}

\end{theorem}

\begin{rmk}
Let $H$ and $A$ be groups definable over $\Q$. When we say that   ``$H$ is \emph{virtually} $A$'', we mean that ``there are a finite index $\Q$-definable subgroup $X$ of $H$ and a $\Q$-definable morphism $f: X\rightarrow A$ such that both $\ker(f)$ and $\img(f)$ are finite''. So the first part of Theorem \ref{sta-no-Ga2} says that $H$ is  virtually a finite index subgroup of a product of $\Ga^s\times\Gm^t$.
\end{rmk}

We also study the $\mu$-invariance of weakly generic types on  $\dfg$ groups over $\Q$. Briefly, letting $\mu$ be the partial type consisting of all definable open neighborhoods of $\text{id}_G$ over $\Q$ and $p\in S_G(M)$ a weakly generic type on $G$, we call $p$  {\em $\mu$-invariant} if  $\mu\cdot p=p$. The  $\mu$-invariance is interesting in its own right, because a weakly generic type on a $\dfg$ group looks like something invariant under ``small'' disturbances. We will show that

\begin{theorem}\label{mu-inv-1}
Let $H$ be a commutative $\dfg$ group definable over $\Q$. Then any weakly generic $p\in S_H(M)$ is almost $\mu$-invariant.
\end{theorem}

Finally, we will give an answer to the Newelski's Question for commutative  groups definable over $\Q$, in both local and global cases. What we call the global context is where types are over a monster model $\M$. What we call the local context is where $M$ is any elementary extension of $\Q$, and we pass to the Shelah expansion $M_0=M^\ext$ of $M$  by adding externally definable sets and consider instead the action of $G(M_0)$ on $S_G(M_0)$.
\begin{theorem}\label{Main}
Let $G$ be a commutative $\Q$-definable group which is not $\dfg$. Then
\begin{enumerate}[(i)]
\item (Local case)  $\AP(S_G(M_0))=\WGen(S_G(M_0))$ if and only if $G$ is  stationary.
\item (Global case) $\AP(S_G(\M))=\WGen(S_G(\M))$ if and only if $G$ is stationary.
\end{enumerate}
\end{theorem}

\subsection{Notations and conventions}

$L$ will denote a language, $T$  a complete theory, $\M$  a monster model, and $M,N,...$ small submodels of $\M$, where we call a subset/submodel $A$ of $\M$  small if $|A|<|\M|$. As the theories which we consider are $\NIP$, we always assume that $T$ is $\NIP$. Every type over a small subset is realized in $\M$.
We call $N$ a {\em sufficiently saturated} submodel of $\M$ if $N$ is $\lambda$-saturated for some sufficiently large cardinal $\lambda$, and, of course, every type over $N$ is realized in $\M$.
We usually write tuples as $a, b, x, y...$ rather than $\bar a, \bar b, \bar x, \bar y...$. Let $\phi(x)$ be any $L_\M$-formula with $x=(x_1,..,x_n)$, and $A\sq \M$, then $\phi(A)$ is defined to be the set $\{a\in A^n|\ \M\models \phi(a)\}$. A ``type'' is a complete type, and a ``partial type'' is a partial type. By a ``global type'' we mean a complete type over $\M$.  Let $A\sq B$ and $p$  a (partial) type over $B$, then $p\upharpoonright_A=\{\phi(x)\in p|\ \phi\in L_A\}$ denotes the restriction of $p$ to $A$. Let $L'\supseteq L$ be an expansion of $L$ and $M_0$ an $L'$-structure, by $M_0\upharpoonright_L$ we means the reduct of $M_0$ to $L$.

When we speak of a set definable {\em in} $M$, $N$ or $\M$, we mean the obvious thing.  When we speak of a set $X$ definable {\em over}
$M$ we typically mean a set definable in $\M$ defined with parameters from $M$. We sometimes use $X(x)$ to denote the formula which defines $X$. In this case, $X(M)$ denotes the set definable in $M$ by the same formulas defining $X$ in $\M$.
On the other hand if for example $X$ is a set definable {\em in} $M$ then we can consider the sets definable in $N$ or $\M$ by the same formulas, which we may refer to as $X(N)$, $X(\M)$.

Let $\Sigma(x)$ be a partial type over a small subset, abusing
notations, we will identify $\Sigma$ with the realizations of $\Sigma(x)$ in $\M$, namely $\Sigma=\Sigma(\M)$. A set $Y$ is type-definable in $\M$ if it is defined by a partial type over a small subset.  Similarly, by a type-definable object, we mean a type-definable object in the monster model $\M$.

Our notations for model theory are standard, and we will assume familiarity with basic notions such as type spaces, heirs, coheirs, definable types etc. References are \cite{Pozit-Book} and \cite{Sim-Book}.

The paper is organized as follows:

For the rest of this section, we give precise definitions and preliminaries relevant to our results.

In section 2, we will prove Theorem \ref{sta-no-Ga2}.

In section 3, we will prove Theorem \ref{mu-inv-1}.

In section 4, we will prove Theorem \ref{Main}, the main theorem of the article.

\vspace{0.3cm}

\begin{acknowledgment}
  The authors were supported by the National Social Science Fund of China (Grant No.\@ 20CZX050).
\end{acknowledgment}

\subsection{Shelah expansion}

By an externally deﬁnable set in $M$ we mean a subset of $M^n$ of the form $\phi(N,b)\cap M^n$ (i.e. trace on $M^n$) where $\phi(x,y)$ is an $L$-formula, $N\succ M$, and $b\in N^k$ with some $n,k\in\N$.  The \emph{Shelah expansion} $M^\ext$ of $M$ is the expansion of $M$ obtained by adding predicates for all externally deﬁnable sets in $M$. We denote the language of $M^\ext$ by $L^\ext_M$ and the theory of $M^\ext$ by $T_M^{\ext}$.  Let $G$ be a group definable over $M$. We denote the collection of quantifier-free types over $M^\ext$ which concentrate on $G$ by $S_{G,\ext}(M)$. The space $S_{G,\ext}(M)$ is naturally homeomorphic to the space $S_{G,M}(\M)$ of global types  concentrating on $G$ which are finitely satisfiable in $M$, via the map
\[
 S_{G,M}(\M)\ni p\mapsto \{\psi(M)|\ \psi\in p \}\in S_{G,\ext}(M).
\]

As $T$ has $\NIP$,  we have:
\begin{fact}\label{fact-SG(M^ext)}\cite{Sim-Book}, Proposition 3.23 and Corollary 3.24]
Assume that $M$ is a model of $T$, then
\begin{enumerate}[(i)]
    \item $T^\ext_M$ has quantifier elimination. So $S_{G,ext}(M)$ coincides with  $S_G(M^\ext)$, the space of types over $M^\ext$ concentrating on $G$.
    \item $T^\ext_M$ also has $\NIP$.
    \item ${(M^\ext)}^\ext=M^\ext$, namely, every type over $M^\ext$ is definable.
\end{enumerate}
\end{fact}

Note that if all types over $M$ are definable (for example, $M=\Q$, see \cite{Delon}), then every externally definable  set in $M$ is already definable, so $M=M^\ext$.

\begin{fact}[Exercise 3.10, \cite{Sim-Book}]\label{Sim-Exercise}
Let $N_0$ be an elementary extension of $M^\ext$ and $N=N_0\rd_L$ the reduct of $N_0$ to $L$. Then there is a set $B\subset\M$ of parameters such that every definable set in $N_0$ is externally definable in $N$ by a formula in $L_B$.
\end{fact}
With the notaions as in the Fact \ref{Sim-Exercise}, we have that
$S_G(N_0)$ is naturally homeomorphic to $S_{G,N}(B)$, the space of types over $B$ concentrating on $G$ which are finitely satisfiable in $N$.

\subsection{Definable topological dynamics}

Our reference for (abstract) topological dynamics are \cite{book-mini-flows-I,book-mini-flows-II}. Let $G$ be a topological  group, by a \emph{$G$-flow} we mean  an  action $G\times X\rightarrow X$ of $G$ on  a  compact  Hausdorff  topological space $X$ by homeomorphisms, and denote it by $(G,X)$. We often assumed that $X$ has a dense $G$-orbit.

A \emph{subflow} of $X$ is a closed $G$-invariant subset of $X$. Minimal(under the relation of inclusion) subflows of $X$ always exist by Zorn's Lemma.
A point $x\in X$ is \emph{almost periodic} if $\cl(G\cdot x)$, the closure  of its $G$-orbit, is  a  minimal subflow of $X$ or, equivalently, if $x$ is in some minimal subflow of $X$. We denote the space of all almost periodic points in $X$ by $\AP(X)$, which is also the union of all minimal subflows of $X$.

Given a $G$-flow $(G,X)$, its \emph{enveloping semigroup} $E(X)$ is the closure in the space $X^X$ (with the product topology) of the set of maps $\pi_g: X\rightarrow X$ for $g\in G$, where $\pi_g(x)=g\cdot x$, equipped with the composition (which is continuous on the left). So any $e\in E(X)$  is
 a map from $X$ to $X$.

\begin{fact}\label{Ellis-env}
Let $X$ be a $G$-flow. Then
\begin{enumerate}[(i)]
    \item $E(X)$ is also a $G$-flow and $E(E(X))\cong E(X)$ as $G$-flows.
    \item For any $x\in X$, the closure of its $G$-orbit is exactly $E(X)(x)$. Particularly, for any $f\in E(X)$, $E(X)\circ f$ is the closure of $G\cdot f$.
\end{enumerate}
\end{fact}

\begin{fact}\label{AP-bi-ideal}
Let $X$ be a $G$-flow. Then $\AP(E(X))$ is a bi-ideal of $E(X)$.
\end{fact}
\begin{proof}
For each $f\in \AP(E(X))$, $E(X)\circ f$ is a minimal subflow, so $E(X)\circ f\sq \AP(E(X))$, and thus $\AP(E(X))$ is a left ideal.

We now show that $\AP(E(X))$ is also a right ideal. Let $f\in \AP(E(X))$ and $g\in E(X)$. Then it suffices to show that $f\circ g$ is almost periodic. It is easy to see that
\[
\cl(G\cdot (f\circ g))=E(X)\circ (f\circ g)=(E(X)\circ f)\circ g.
\]
For any $h\in (E(X)\circ f)$, we have that $E(X)\circ h=E(X)\circ f$ since $E(X)\circ f$ is minimal. So
\[
\cl(G\cdot (h \circ g))=E(X)\circ (h\circ g)=(E(X)\circ h)\circ g=(E(X)\circ f)\circ g=\cl(G\cdot (f\circ g)).
\]
We conclude that $\cl(G\cdot (f\circ g))$ is minimal since $\cl(G\cdot h')=\cl(G\cdot (f\circ g))$ for each $h'\in \cl(G\cdot (f\circ g))$, and hence $f\circ g$ is almost periodic as required.
\end{proof}


In the model theoretic context,  we consider a  group $G$ definable over $M$ and the action of $G(M)$ on its type space $S_G(M)$ as $gp=\tp(ga/M)$ where $g\in G(M)$ and  $a$ realizes $p$. It is easy to see that $S_G(M)$ is a $G(M)$-flow with a dense orbit $\{\tp(g/M)|\ g\in G(M)\}$.

Take a monster model $\M$ and identify $G$ with $G(\M)$. We call a formula $\varphi(x)$, with parameters in $\M$, a \emph{$G$-formula} if $\varphi(\M)$ is a definable subset of $G$. A partial type $r$ is called a {\em $G$-type} if every formula in $r$ is a $G$-formula.  Suppose that $\varphi(x)$ is a $G$-formula  and $g\in G$, then the left translate $g\varphi(x)$ is defined to be $\varphi(g^{-1}x)$.  It is easy to check that $(g\varphi)(\M)=g\varphi(\M)$.  For $p\in S_G(M)$, we have $gp=\{g\varphi(x)|\ \varphi\in p\}$.

We recall some notions from \cite{Newelski-I}.
\begin{definition}
\begin{enumerate}
    \item A definable subset $X\sq G$ is (left) generic if finitely many left translates of $X$ cover $G$. Namely, there are $g_1,...,g_n\in G$ such that $\bigcup_{i=1}^n g_iX=G$.
    \item A definable subset $X\sq G$
    is weakly  generic if there is a non-generic definable subset $Y$ such that $X\cup Y$ is generic.
    \item A $G$-formula $\varphi(x)$ is generic(weakly generic) if the definable set $\varphi(\M)$ is generic(resp. weakly generic).
    \item A (partial) $G$-type  $p$ is generic(weakly generic) if every formula in $p$ is generic(resp. weakly generic).
\end{enumerate}
\end{definition}

Let $\AP(S_G(M))\sq S_G(M)$ be the space of almost periodic types, $\WG(S_G(M))\sq S_G(M)$ the space of weakly generic types, and $\Gen(S_G(M))\sq S_G(M)$ the space of generic types. Then $\WG(S_G(M))=\cl(\AP(S_G(M)))$. Moreover, if $\Gen(S_G(M))\neq \emptyset$, then $\WG(S_G(M))=\Gen(S_G(M))$ is the unique minimal subflow of $S_G(M)$  (see \cite{Newelski-I}).

\begin{fact}\label{Semigp-struc}\cite{Newelski-I}
The enveloping semigroup $E(S_{G, \ext}(M)) $ of  $S_{G, \ext}(M)$ is  isomorphic to $(S_{G,M}(\M),*)$ where $*$ is defined as following: for any $p,q\in S_{G,M}(\M)$, $p *q=tp(b\cdot c/\M)$ with $b$ realizes $p$ and $c$ realizes $q$, and $tp(b/\M,c)$ is finitely satisfiable in $M$.
\end{fact}

\begin{rmk}
Note that $S_{G,M}(\M)$ is naturally homeomorphic to $S_{G,\ext}(M)$.  Assume that $T$ has NIP, we see from Fact \ref{fact-SG(M^ext)} and Fact \ref{Semigp-struc} that there is a semigroup operation ``$*$'' on $S_G(M^\ext)$, which is  defined as following: for any $p,q\in S_G(M^\ext)$, $p *q=tp(b\cdot c/M^\ext)$ with $b$ realizes $p$ and $c$ realizes the unique heir of  $q$ over $ \dcl(M^\ext,b)$. Moreover, it is easy to see from Fact \ref{Ellis-env}  that for any $p\in S_{G}(M^\ext)$,
\[
\cl(G(M)\cdot p)= S_{G}(M^\ext) *p.
\]
\end{rmk}

\subsection{$\NIP$, definable amenability, and connected components}

 Recall that $T$ has $\NIP$,  which means that,  for any indiscernible sequence $(b_i:\ i<\omega)$, formula $\phi(x,y)$, and $a\in \M$, there is an eventual truth-value of $\phi(a,b_i)$ as $i\rightarrow \infty$.

Let $G=G(\M)$ be a definable group. Recall that a type-definable over $A$ subgroup $H$ is  a  type-definable over $A$ subset of $G$, which is also a subgroup of $G$. We say that $H$ has bounded index if $|G/H|<2^{|T|+|A|}$. For groups definable in  $\NIP$ structures, the smallest type-definable subgroup of bounded index exists \cite{HPP-$NIP$-JAMS}, which is the intersection of all type-definable subgroups of bounded index, we write it as $G^{00}$, and call it the \emph{type-definable connected component}. Another model theoretic invariant is $G^{0}$, called the \emph{definable-connected component} of $G$, which is the intersection of all definable subgroups of $G$ of   finite index.  Clearly, $G^{00}\leq G^0$.

\begin{fact}\cite{CPS-Ext}\label{G00=G00ext} If $M\models T$ and $G$ is definable over $M$,  then $G^{00}$ is the same whether computed in $T$ or $T^{\ext}_M$.
\end{fact}

\begin{fact}\label{times00}
Let $G_1$ and $G_2$ be definable groups and $G=G_1\times G_2$. Then $G^{00}=G_1^{00}\times G_2^{00}$. Here, we identify $G_1$ with $G_1\times \{1_{G_2}\}$ and $G_2$ with $\{1_{G_1}\}\times G_2$, subgroups of $G$.
\end{fact}
\begin{proof}
Because $[G_1\times G_2:G_1^{00}\times G_2^{00}]=[G_1:G_1^{00}][G_2:G_2^{00}]$ is bounded, we have that $G^{00}\leq G_1^{00}\times G_2^{00}$. For $i=1,2$, since ${G_i}/{(G_i\cap G^{00})}\cong {(G^{00}G_i)}/{G^{00}} \leq {G}/{G^{00}}$, $G_i\cap G^{00}$ is a type definable subgroup of $G_i$ with bounded index,  we see that $G_i^{00} \leq G_i\cap G^{00}$, and thus $G_1^{00}\times G_2^{00}\leq G^{00}$.
\end{proof}

Recall also that a  \emph{Keisler measure} over $M$ on $X$, with $X$ a definable subset of $M^n$, is a finitely additive measure on the Boolean algebra of $M$-definable subsets of $X$. When we take the monster model, i.e. , $M=\M$,  we call it a global Keisler measure. A definable group $G$ is said to be \emph{definably amenable} if it admits a global (left) $G$-invariant probability  Keisler measure. By \cite{HPP-$NIP$-JAMS} this is equivalent to the existence of a $G(M)$-invariant probability  Keisler measure over $M$ on $G$, whenever $M$ is a model over which $G$ is defined.
Note that if $G$ is commutative, then $G$ is amenable, so is definably amenable.

\subsection{
Strongly $f$-generics, definable $f$-generics, and finitely satisfiable generics}

Let $\phi(x,y)$ be a formula. Recall that a formula $\phi(x,b)$ \emph{divides} over a set $A$ if there is an infinite $A$-indiscernible sequence $(b=b_0,b_1,b_2,...)$ such that $\{\phi(x,b_i)|\ i<\omega\}$ is inconsistent. A type $p\in S(B)$ divides over $A\sq B$ if there is a formula $\phi\in p$ divides over $A$.  Since $T$ has  $\NIP$, a global type $p\in S(\M)$ does not divide over   $M$ if and only if $p$ is $\Aut(\M/M)$-invariant (see \cite{HP-$NIP$-inv-measure}).

Let $G=G(\M)$ be a definable group and $p\in S_G(\M)$ a global type on $G$. Following the notation of \cite{CS-amenable-$NIP$-group}, we call a global type $p\in S_G(\M)$   \emph{strongly $f$-generic over $M$}  if  $gp$ does not divide over $M$ for each $g\in G$. We call   $p\in S_G(\M)$ a strongly $f$-generic type if it is strongly $f$-generic over some small submodel.
A nice result of \cite{CS-amenable-$NIP$-group} shows that:

\begin{fact}\label{G00-inv=WG} Let $G$ be a definably amenable group and $p$ be a global $G$-type. Then
\begin{enumerate}[(i)]
    \item $p$ is weakly generic iff it is $G^{00}$-invariant iff it has a bounded $G$-orbit.
    \item $p$ is strongly $f$-generic over $M$ iff it is weakly generic and $M$-invariant (or does not divide over $M$).
\end{enumerate}
\end{fact}

Among the strongly $f$-generics $p\in S_G(\M)$, there are two extreme cases:
\begin{enumerate}
\item There is a small submodel $M$ such that every left $G$-translate of $p$ is  finitely satisfiable in $M$, and we call such types the \emph{$\fsg$} (finitely satisfiable generic) types on $G$ over $M$;
\item There is a small submodel $M$ such that every left $G$-translate of $p$ is definable over $M$, and we  call such types the \emph{$\dfg$} (definable $f$-generic) types on $G$ over $M$.
\end{enumerate}

A definable group $G$ is called $\fsg$ or $\dfg$ if is has a global $\fsg$ or $\dfg$ type, respectively.
Both $\fsg$ and $\dfg$ groups are definably amenable(see \cite{CS-amenable-$NIP$-group}).
We now discuss these two cases.

By \cite{HPP-$NIP$-JAMS} we have:
\begin{fact}\label{fact-fsg}
 Let $G$ be an $\emptyset$-definable $\fsg$ group witnessed by a $\fsg$ type $p\in S_G(\M)$. Then
\begin{enumerate}[(i)]
    \item $p$ is both left and right generic. So the space $\Gen(S_G(\M))$ of global generic types is the unique minimal subflow of $S_G(\M)$.
    \item Any left (right) translate of $q\in \Gen(S_G(\M))$ is finitely satisfiable in any small submodel of $\M$.
    \item For any $N\prec\M$, every generic type $q\in S_G(N)$ has a unique global generic extension.
    \end{enumerate}
    \end{fact}

\begin{fact}\cite{CPS-Ext}\label{fact-fsg-ext}
 Let  $C$ be an $\emptyset$-definable $\fsg$ group and $M\prec\M$.
    \begin{enumerate}[(i)]
    \item $C$ also has $\fsg$ when we compute it in $T^\ext_M$.
    \item  $q\mapsto \{\psi(M)|\ \psi\in q\}$ is a bijection between $\Gen(S_C(\M))$ and $\Gen(S_C(M^\ext))$.
\end{enumerate}
\end{fact}

We now discuss the $\dfg$ groups.

\begin{fact}\label{dfg-fact-1}\cite{Pillay-Yao-mini-flow}
Let   $G$ be an $\emptyset$-definable group and $p\in S_G(\M)$ a  weakly generic type. If $p$ is definable over  $M$, then
\begin{enumerate}[(i)]
    \item Every left  translate of $p$ is definable over $M$;
    \item $G^{00}=G^{0}$.
    \item $G\cdot p$ is closed, and hence a minimal subflow of $S_G(\M)$.
\end{enumerate}
\end{fact}

\begin{fact}\label{dfg-fact-2}\cite{CPS-Ext}
 Let $G$ be a $\dfg$ group definable over $M$. Then $G$ also has $\dfg$ when we compute it in $T^\ext_M$.
\end{fact}

\begin{fact}\label{min-flow-restriction}\cite{Pillay-Yao-mini-flow}
Assume that $G$ is a definably amenable group definable over $M$ and $M\prec N$.
\begin{enumerate}
    \item [(i)] Let $\pi: S_G(N^\ext)\rightarrow S_G(N)$ be the canonical restriction map, and ${\cal M}$   a minimal $G(N)$-subflow of $S_G(N^\ext)$. Then $\pi({\cal M})$ is a  minimal $G(N)$-subflow of $S_G(N)$.
    \item [(ii)]  Let $\pi: S_G(N)\rightarrow S_G(M)$ be the canonical restriction map, and $\cal M$  a minimal $G(N)$-subflow of $S_G(N)$. Then $\pi({\cal M})$ is a  minimal $G(M)$-subflow of $S_G(M)$.
\end{enumerate}
\end{fact}

\subsection{Commutative Definable groups over $\Q$}

Let $p$ be a prime and $\Q$ the field of $p$-adic numbers. We call the complete theory of $\Q$, in the language of rings, the theory of \emph{$p$-adically closed fields}, written $\pCF$. Macintyre showed in \cite{Mac} that $\pCF$ has quantifier elimination in the language of rings together with new predicates \emph{$P_n(x)$} for the $n$-th powers for each $n\in \N^{>0}$. It is well known that $\pCF$ has $\NIP$(see \cite{Sim-Book}). As  $\pCF$ has definable Skolem functions \cite{Dries-skolem-functions}, every definable map has a definable section. The valuation group $(\mathbb Z, +,<)$ and the valuation $v : \Q \rightarrow \mathbb Z\cup\{\infty\}$ are interpretable. A $p$-adically closed field is a model of $\pCF$. For any $M\models \pCF$, 
$\Gamma_M$ will denote the value group, which is an elementary extension of $(\mathbb{Z},+,<)$. If $A\sq \Gamma_M$, we call $\alpha\in\Gamma_M $ is \emph{unbounded over $A$} if either $\alpha<\dcl(A)$ or $\alpha>\dcl(A)$. The relation $v(x)\leq v(y)$ is $\emptyset$-definable in the language of rings (see \cite{Denef-I}). The $p$-adic field $\Q$ is a locally compact topological field, with basis given by the
sets
\[
B(a, n) = \{x \in \Q |  v(x-a)\geq n\}
\]
for $a\in\Q$ and $n\in\mathbb Z$. We now assume that $T=\pCF$ is the theory of $\Q$ in the language of rings and $\M$ a monster model of $\pCF$.

 An \emph{definable $p$-adic analytic manifold} is a topological space with a covering by finitely many open sets each homeomorphic to an open definable (in $\Q$) subset of $\Q^{n}$ for some $n$ such that the transition maps  are definable and analytic.   A \emph{definable $p$-adic analytic group} is a definable $p$-adic analytic manifold equipped with a group structure which is definable and analytic when read in the appropriate charts. Such a definable $p$-adic analytic group is a definable group in $\Q$. Conversely, by Lemma 3.8 of \cite{Pillay-On fields definable in Qp}, any group $G$ definable in $\Q$ can be definably equipped with the structure of a definable $p$-adic analytic group. As indicated in \cite{O-P}, that any $p$-adic analytic group $G$ has a family of definable open compact subgroups which forms a local base of the identity $\id_G$ of $G$.

Recall from \cite{J-Y-On noncompact p-adic groups} that a definable $p$-adic analytic manifold $X$ (in the monster model) is {\em definably compact} if the following holds:  for any definable family ${\cal F}=\{Y_t|\  t\in T\}$ of non-empty closed sets $Y_t \sq X$, if $\cal F$ is downwards directed, then $\bigcap {\cal F}\neq \emptyset$.
Note that $X$ is definably compact iff $X(\Q)$ is compact when $X$ is definable over $\Q$.
If $X$ is definable over $\Q$ and definably compact, then  for any $a^*\in X$  there is a unique $a\in X(\Q)$ such that $a^*\in U$  for every $\Q$-definable open subset $U\subset X$ containing $a$. We write $a=\st(a^*)$, and call  $\st: X\rightarrow X(\Q)$ the standard part map.

\begin{fact}\cite{JY-abelian-pcf-group}\label{dfg-fsg-decomp}
Let $G$ be an abelian group definable over $\Q$. Then $G^{00}=G^0$ and there is a definable short exact sequence
\begin{equation}\label{dfg-fsg-short-ext-seq}
    1\rightarrow H\rightarrow G\rightarrow_\pi C\rightarrow 1
\end{equation}
with $H$ a $\Q$-definable $\dfg$ group,  $C$ a $\Q$-definable $\fsg$ group, and $\pi$ a $\Q$-definable homomorphism.
\end{fact}

\begin{fact}\label{fsg-G00=G0}\label{fsg=compact}
Let $C$ be a group definable over $\Q$.
\begin{enumerate}
  \item [(i)] $C$ is definably compact iff $C$ has $\fsg$.
  \item  [(ii)] If $C$ is definably compact, then $C^{00}=C^0$ coincides with $\ker(\st)$, and $ \st$ induces a homeomorphism between $C/C^{0}$ (with its logic topology) and the $p$-adic analytic group $C(\Q)$.
\end{enumerate}
\end{fact}
\begin{proof} For (i) the right implies left direction of (ii) is Corollary 2.3 (iv) of \cite{O-P}. The left to right  appears in  Proposition 3.1 of \cite{Johnson-fsg}. For (ii) see Corollary 2.4 of \cite{O-P}.
\end{proof}


Let $F$ be an algebraic closed field containing $\Q$. We call $G(F)$ an algebraic group over $\Q$ if the variety structure as well as the group structure are given by data (polynomial equations, transition maps, morphisms over $\Q$ (See \cite{Pillay-ACF}). A  {\em $p$-adic algebraic group} is the group $G(\Q)$ of $\Q$-points of an algebraic group $G(F)$ over $\Q$.  Of course, $G(\Q)$ will be also a definable group  in the structure  $(\Q,+,\times,0,1)$,
but essentially just quantifier-free definable in the ring language.  Abusing terminology, we also refer to $G(\M)$ as a $p$-adic algebraic group if $G(\Q)$ is. To consistent with our earlier notation, we denote $G(\M)$ by $G$. We denote the additive group $(\M,+)$ by $\Ga$ and the multiplicative group  $(\M^*,\times)$ by $\Gm$. We denote the direct product of $n$ copies of $\Ga$ (resp. $\Gm$) by $\Ga^n$ (resp. $\Gm^n$).
We call a (connected) $p$-adic algebraic group $G$ {\em trigonalizable} over $\Q$ if there is a normal sequence
\[
\{1_G\}=G_0\vartriangleleft ... \vartriangleleft G_i \vartriangleleft G_{i+1} \vartriangleleft... \vartriangleleft G_n=G
\]
such that each $G_{i+1}/G_i$ is (quantifier free) definably isomorphic to  $\Ga$ or $\Gm$ over $\Q$.

It was showed in \cite{PY-dfg-groups} that every $\dfg$ group over $\Q$ is virtually a finite index subgroup  of a $p$-adic algebraic groups trigonalizable  over $\Q$, precisely, we have

\begin{fact}\label{Almost-alge}\cite{PY-dfg-groups}
Let $H$ be a $\dfg$ group definable over $\Q$. Then $H$ is virtually  a connected algebraic group  which is trigonalizable over $\Q$.
\end{fact}

For commutative  algebraic groups, we have

\begin{fact}[Theorems 17.17 and Corollary 17.19, \cite{J.milne}]\label{Springer}
Let $G$ be a connected commutative linear algebraic group over $\Q$. Then $G(\Q)$ is definably isomorphic to $G_u(\Q)\times T(\Q)$ over $\Q$, where $T(\Q)$ is a $\Q$-tours and $G_u(\Q)$ is a product of copies of $\Ga(\Q)$.
\end{fact}

\begin{rmk}\label{split-tours}
 Any $\Q$-tours $T(\Q)$ is an almost direct product of $T_{spl}(\Q)$ and $T_{an}(\Q)$, where $T_{spl}(\Q)$ is $\Q$-split (i.e. $\Q$-definably isomorphic to a product of copies of $\Gm(\Q)$), and $T_{an}$ is anisotropic (see \cite{V. Platonov-A. Rapinchuk}, p. 53), and thus is compact (see \cite{SP. Wang}).
\end{rmk}

\begin{fact}\cite{P-Y-commu-by-fin}\label{commu-open-neibor}
Let $G$ be a group definable in $\Q$. Suppose that $G$ has a commutative open neighborhood of $1_G$, then $G$ is commutative-by-finite.
\end{fact}

\begin{coro}\label{H-product}
Let $H$ be a commutative $\dfg$ group definable over $\Q$. Then $H$ is virtually  a product of copies of $\Ga$ and $\Gm$.
\end{coro}
\begin{proof}
By Fact \ref{Almost-alge}, there are a finite index $\Q$-definable subgroup $A<H$ and a finite subgroup $A_0<A$ such that $A/A_0$ is isomorphic to an open subgroup of $B$ where $B$ is a $p$-adic algebraic group trigonalizable over $\Q$. By Fact \ref{commu-open-neibor}, $B$ is commutative-by-finite, so is commutative as $B$ is connected.

Now $B$ is a commutative $p$-adic linear algebraic group, by Fact \ref{Springer}, $B(\Q)$ is  definably isomorphic to $\Ga^s(\Q)\times T(\Q)$ over $\Q$, where $s\in \N$ and $T(\Q)$ is a $\Q$-tours. Since $B$ is trigonalizable over $\Q$, we see from Remark \ref{split-tours} that $T(\Q)$ is split over $\Q$, thus is definably isomorphic to $\Gm^t(\Q)$ over $\Q$ for some $t\in \N$. This completes the proof.
\end{proof}

\section{Stationarity and badness}

Recall that we work in $\pCF$.  We introduce the stationarity first.

\begin{definition}\label{def-of-stationary}
We say that a $\Q$-definable group $G$ is {\em stationary} if every
weakly generic type  $p\in S_G(\Q)$  has just one global weakly generic extension.
\end{definition}

In \cite{YZ-Sta}, we consider the case where $G$ is a $\Q$-definable group (not necessarily commutative) admitting a $\Q$-definable short exact sequence
\[
1\rightarrow H\rightarrow G\rightarrow_\pi C\rightarrow 1,
\]
with $C$ a $\fsg$ group and $H$ a $\dfg$ group, we proved in  \cite{YZ-Sta} that:

\begin{fact}\label{sta-AP=WG}

\begin{enumerate}
\item [(i)] $G$ is stationary iff $H$ is stationary.
\item [(ii)] (Local case)  If $G$ be stationary, then $\AP(S_G(M^\ext))=\WGen(S_G(M^\ext))$ for any model $M\succ \Q$.
\item [(iii)] (Global case)  If $G$ be stationary, then $\AP(S_G(\M))=\WGen(S_G(\M))$.
\end{enumerate}
\end{fact}

For a $\dfg$ group definable over $\Q$, we have that
\begin{fact}\cite{YZ-Sta}\label{H-sta-equi}
Let $H$ be a dfg group definable over $\Q$, then the following are equivalent:
\begin{enumerate}
    \item [(i)] $H$ is stationary.
    \item  [(ii)] $H$ has boundedly many global weakly generic types.
    \item [(iii)]  There is a small model $M$ such that every global weakly generic type is $M$-definable.
    \item [(iv)]  Every global weakly generic type is $\Q$-definable.
\end{enumerate}
\end{fact}

Now we are going to study the stationarity of commutative groups definable over $\Q$. Firstly, we show that the stationarity is preserved by the relation of ``virtually'' on definably amenable groups.

Firstly, the relation of being a finite index definable subgroup does not matter. This because, for definable groups $H<G$ with $[G:H]$ finite, we have that $H^{00}=G^{00}$ and any weakly generic type on $G$ is indeed on $H$.

Now we only need to deal with the case when $A$ is a quotient of $G$ with a finite kernel.

\begin{lemma}\label{finite-to-one}
    Let $G$ and $A$ be groups definable over a model $M$ and $\pi:G\rightarrow A$  a  surjective morphism definable over $M$.  If $p\in S_A(M)$ and $q=\pi(p)$, then $\pi^{-1}(q) =\ker(\pi)p$.
\end{lemma}
\begin{proof}
Let $r\in \pi^{-1}(q)$. Suppose that $a\models r$, then $\pi(a)\models q$. Let $f$ be a definable section of $\pi$, then we have that $a\in \ker(\pi)f(\pi(a))$, which means that $r\in \ker(\pi)f(q)$, so $\pi^{-1}(q)\subseteq\ker(\pi)f(q)$. On the other side, it is easy to see that  $\ker(\pi)f(q)\sq \pi^{-1}(q)$, so $\ker(\pi)f(q)= \pi^{-1}(q)$. As $p\in \pi^{-1}(q)$, we have $\ker(\pi)p=\ker(\pi)f(q)=\pi^{-1}(q)$.
\end{proof}

\begin{fact}\cite{J-Y abelian groups}
  Let $\pi : G \to A$ be a definable surjective morphism of definable groups.  Then
  $\pi(G^{00}) = A^{00}$.
\end{fact}

\begin{lemma}\label{finite-to-one-G00(N)-inv}
Let $\pi: G\rightarrow A$ be a $\Q$-definable surjective morphism of $\Q$-definable groups. Suppose that $N$ is a sufficiently saturated small submodel of $\M$ and $p\in S_G(\M)$. Then we have
\begin{enumerate}
    \item [(i)] If $p$ is $G^{00}(N)$-invariant, then $\pi(p)$ is $A^{00}(N)$-invariant.
    \item [(ii)] If $\pi(p)$ is $A^{00}(N)$-invariant, then $G^{00}(N)p\sq\ker(\pi)p$.
\end{enumerate}
\end{lemma}
\begin{proof}
If $p\in S_G(\M)$ is $G^{00}(N)$-invariant, then $\pi(p)$ is $\pi(G^{00})(N)$-invariant, thus is $A^{00}(N)$-invariant.

On the other side, suppose that $\pi(p)=q\in S_{A}(\M)$ is $A^{00}(N)$-invariant, then $\pi^{-1}(q)$ is $\pi^{-1}(A^{00}(N))$-invariant.  Since  $G^{00}(N)\sq \pi^{-1}(A^{00}(N))$, we see that $\pi^{-1}(q)$ is $G^{00}(N)$-invariant. By Lemma \ref{finite-to-one}, $\pi^{-1}(q)=\ker(\pi)p$, so $G^{00}(N)p\sq \ker(\pi)p$.
\end{proof}

\begin{lemma}\label{equiv-stationarity-I}
Let $G$ and $A$ be definably amenable groups definable over $\Q$ and $\pi: G\rightarrow A$ a $\Q$-definable surjective morphism with a finite kernel. Let $M\succ \Q$. Then $p\in S_G(M)$ is weakly generic iff $\pi(p)\in S_A(M)$ is weakly generic.
\end{lemma}
\begin{proof}
Since every weakly generic type over $M$ extends to a global weakly generic type, it suffices to show that $p\in S_G(\M)$ is weakly generic iff $\pi(p)\in S_A(\M)$ is weakly generic.

If $p\in S_G(\M)$ is weakly generic, then $p$ is $G^{00}$-invariant, so $\pi(p)$ is $A^{00}$-invariant by part (i) of Lemma \ref{finite-to-one-G00(N)-inv}. On the other side, suppose that  $q=\pi(p)\in S_{A}(\M)$ is weakly generic, where $p\in S_G(\M)$. Then $q$ is $A^{00}$-invariant, and by part (ii) of Lemma \ref{finite-to-one-G00(N)-inv} $G^{00}p\sq \ker(\pi)p$ is finite. So $p\in S_G(\M)$ has a bounded $G$-orbit and thus is weakly generic  by Fact \ref{G00-inv=WG}.
\end{proof}

\begin{rmk}\label{rmk-ker(f)}
    Suppose that $G$ is virtually a finite index subgroup of $A$, witnessed by a $\Q$-definable finite index subgroup $X\sq G$ and  a $\Q$-definable morphism $f: X\rightarrow A$. We will always assume that $\ker(f)\sq G^0(\Q)$.
\end{rmk}
\noindent
{\em Explanation.} Since $\ker(f)$ is finite and $G^0$ is the intersection of all finite-index subgroup definable over $\Q$, there is  a $\Q$-definable finite index subgroup $Y$ of $X$ such that $Y\cap \ker(f)=G^0\cap \ker(f)$.  We see that the kernel of $f\upharpoonright_Y:Y\rightarrow A$ is contained in $G^{0}$. Replacing $X$ by $Y$ if necessary, we may assume that $\ker(f)\sq G^{0}\cap G(\Q)=G^{0}(\Q)$.

\begin{lemma}\label{equiv-stationarity-II}
Let $G$ and $A$ be definably amenable groups, both definable over $\Q$. If $G$ is virtually a finite index subgroup of $A$ and $G^0=G^{00}$, then $G$ is stationary iff $A$ is stationary.
\end{lemma}
\begin{proof}
Let $X\sq G$ be a $\Q$-definable finite index subgroup and $f: X\rightarrow A$ a definable morphism  such that $\img(f)$ has finite index in $A$.
It is easy to see that $G$ is stationary iff $X$ is stationary, and $A$ is stationary iff $\img (f)$ is stationary. So it suffices to show that $X$ is stationary iff $\img (f)$ is.

If $p\in S_G(M)$ is weakly generic, then $\ker(f)p=p$ by Remark \ref{rmk-ker(f)}. We see from Lemma \ref{finite-to-one} and Lemma \ref{equiv-stationarity-I} that the map
$p\mapsto f(p)$ is a one-one correspondence between $\WGen( S_X(M))$ and $\WGen( S_{\img (f)}(M))$ for arbitrary $M\succ \Q$, which implies that $X$ is stationary iff $\img (f)$ is stationary. This completes the proof.
\end{proof}

Now let $H$ be a commutative dfg group over $\Q$, then it is virtually a finite index subgroup of $\Ga^s\times\Gm^t$ for some $s,t\in\mathbb{N}$. By Lemma \ref{equiv-stationarity-II}, to see the stationarity of $H$, it is reasonable to assume that $H=\Ga^s\times\Gm^t$ since $H^0=H^{00}$.

\begin{fact}\cite{PPY-SL2Qp}\label{Ga0+Gm0}
Let $H$ be either $\Ga$ or $\Gm$. Let $\Sigma_0(x)$ be the partial type $\{v(x)>\gamma|\ \gamma\in \Gamma_\M\}$ and $\Sigma_\infty(x)$ be the partial type $\{v(x)<\gamma|\ \gamma\in \Gamma_\M\}$, then
\begin{enumerate}
    \item [(i)] Every global weakly generic type on $H$ is $\emptyset$-definable.
    \item [(ii)]  If $H=\Ga$, then $H^{00}=H^{0}=H$, and $p\in S_H(\M)$ is weakly generic iff it is consistent with the partial type $\Sigma_\infty(x)$;
    \item  [(iii)] If $H=\Gm$, then $H^{00}=H^{0}=\bigcap_{n\in \N^+}P_n(\Gm)$, and $p\in S_H(\M)$ is weakly generic iff it is either consistent with the partial type $\Sigma_\infty(x)$, or consistent with the partial type $\Sigma_0(x)$.
\end{enumerate}
\end{fact}

We see directly from Fact \ref{Ga0+Gm0} that any global weakly generic type on  $\Ga$  is also a global weakly generic type on $\Gm$. Moreover, we have

\begin{fact}\cite{Yao-Presb}\label{gen-Gm-defble}
Any global weakly generic type on $\Gm^t$ is $\emptyset$-definable for each $t\in\N^+$. Moreover, a global type $\tp(a_1,...,a_t/\M)$  on  $\Gm^t$ is a weakly generic iff $k_1v(a_1)+...+k_tv(a_t)$ is unbounded over $\Gamma_\M$ for all $k_1,...,k_t\in \mathbb Z$ which are not all zero.
\end{fact}

\begin{lemma}\label{gen-Gm-n}
$\tp(a_1,\dots,a_t/\M)$ is a weakly generic type on $\Gm^t$ if and only if each $\tp(a_j/\M,a_i:i<j)$ is $\Gm^{0}$-invariant ($\Gm^0=\Gm^0(\M)$).
Hence, $\Gm^t$ is stationary for each $t\in \N^{+}$.
\end{lemma}
\begin{proof}
Suppose that $\tp(a_1,\dots,a_t/\M)$ is a weakly generic, then it is $(\Gm^0)^t$-invariant, so each $\tp(a_j/\M,a_i:i<j)$ is $\Gm^{0}$-invariant.

we show the other direction by induction on $t\in \N^{+}$. Assume that this lemma holds for $t-1$. Then we have that $\tp(a_1,...,a_{t-1}/\M)$  is a weakly generic type on $\Gm^{t-1}$ by induction hypothesis. Assume for a contradiction that $\tp(a_1,...,a_{t}/\M)$ is not weakly generic, then there are $k_1,...,k_t\in \mathbb Z$ which are not all zero such that $k_1v(a_1)+...+k_tv(a_t)$ is bounded over $\Gamma_\M$. Take $\alpha<\beta\in \Gamma_\M$ such that $\alpha<k_1v(a_1)+...+k_tv(a_t)<\beta$. If $k_t=0$, then $\alpha<k_1v(a_1)+...+k_{t-1}v(a_{t-1})<\beta$, which contradicts to the weak generality of $\tp(a_1,...,a_{t-1}/\M)$.
If $k_t\neq 0$, then take any $b\in \Gm^0$ such that $k_tv(b)<\alpha-\beta$, we have that
\[
k_1v(a_1)+...+k_tv(ba_t)=k_1v(a_1)+...+k_tv(a_t)+k_tv(b)<\beta+\alpha-\beta=\alpha,
\]
which is also a contradiction as $\tp(a_t/\M,a_1,...,a_{t-1})$ is $\Gm^0$-invariant.
\end{proof}

\begin{lemma}\label{GatimesGm^v-sta}
$H=\Ga\times \Gm^t$ is stationary for each $t\in \N$.
\end{lemma}
\begin{proof}
Let $p=\tp(a_0,a_1...,a_t/\M)\in S_H(\M)$ is be weakly generic type.  Then it is easy to see that  $\tp(a_j/\M,a_i:i<j)$ is $\Gm^{0}$-invariant for $j=0,...,t$ since $H^0=\Ga\times (\Gm^0)^t$. By Lemma \ref{gen-Gm-n}, we see that $p$ is a weakly generic type on $\Gm^{t+1}$, so is $\emptyset$-definable by Fact \ref{gen-Gm-defble}. Since every global  weakly generic type on $H$ is $\emptyset$-definable, $H$ is stationary by Fact \ref{H-sta-equi}.
\end{proof}

We recall the notion of badness from \cite{YZ-Sta}:
\begin{definition}\label{bad}
Let $H$ be a $\dfg$ group definable over $\Q$. We say that $H$ is \emph{bad} if there is a strongly $f$-generic type $\tp(a/\M)$ over $\Q$, and a $\Q$-definable function $\theta$ such that $\tp(\theta(a)/\Q)$ is non-algebraic and $\tp(\theta(a)/\M)$ is finitely satisfiable in $\Q$.
\end{definition}

\begin{rmk}\label{bad-is-stationary}
As pointed in \cite{YZ-Sta}, a bad dfg group is not stationary since each non-algebraic global type can not be both definable over $M$ and finite satisfiable in $M$, due to the distality of $\pCF$ (see \cite{Simon-distal} for the details of the distality).
\end{rmk}

\begin{lemma}
    Let $H$ be a $\dfg$ group definable over $\Q$. If $H$ is virtually  $A$, then $H$ is bad iff $A$ is bad.
\end{lemma}
\begin{proof}
Let $X\sq H$ be a $\Q$-definable finite index subgroup and $f:X\rightarrow A$ a $\Q$-definable morphism with both $\ker(f)$ and $A/\img(f)$ are finite.  Applying definable Skolem functions, let $g:\img(f)\rightarrow H$ be a $\Q$-definable section of $f$. If $p\in S_H(\M)$ is a strongly $f$-generic type over $\Q$, then by Lemma \ref{equiv-stationarity-I}, $f(p)$ is weakly generic. Since $p$ is $\Q$-invariant, we have that $f(p)$ is also $\Q$-invariant, and thus is also strongly $f$-generic over $\Q$ by Fact \ref{G00-inv=WG}. Similarly, we can show that if $q\in S_{\img(f)}(\M)$ is strongly $f$-generic over $\Q$, so is $g(q)\in S_H(\M)$. So $X$ is bad iff $\img(f)$ is bad. Clearly, $H$ is bad iff $X$ is, and $A$ is bad iff $\img(f)$ is. This completes the proof.
\end{proof}


In \cite{YZ-Sta}, we showed that $\Ga^2$ is a bad dfg group. It implies that $\Ga^s$ is bad for all $s\geq 2$. Moreover, we have:

\begin{coro}\label{u-geq-2}
    $\Ga^s\times \Gm^t$ is bad and hence not stationary for each $t$ when $s\geq 2$.
\end{coro}
\begin{proof}
    Let $\tp(a/\M)$ be a strongly $f$-generic type of $\Ga^s$ over $\Q$ and $\theta$ a $\Q$-definable function such that $\tp(\theta(a)/\Q)$ is non-algebraic and $\tp(\theta(a)/\M)$ is finitely satisfiable in $\Q$. Let $\tp(b/\M,a)$ be a strongly $f$-generic type on $\Gm^t$ over $\Q$. Then for any $a_0\in \Ga^s$ and $b_0\in{(\Gm^0)}^t$, we see that $\tp(b_0b/\M,a)=\tp(b_0b/\M,a_0+a)$ does not divide over $\Q$ and $\tp(a_0+a/\M)$ does not divide over $\Q$, so $\tp(a_0+a,b_0b/\M)$ does not divide over $\Q$, and thus $\tp(a,b/\M)$ is a strongly $f$-generic type of $\Ga^s\times \Gm^t$ over $\Q$.
    Let $\theta^*(a,b)=\theta(a)$, then $\tp(\theta^*(a,b)/\Q)$ is non-algebraic and $\tp(\theta^*(a,b)/\M)$ is finitely satisfiable in $\Q$.
\end{proof}

Then combine with Lemma \ref{GatimesGm^v-sta}, we have that
\begin{theorem}\label{Sta-s-leq-1}
Let $H$ be an abelian $\dfg$ group definable over $\Q$ which is virtually a finite index subgroup of $\Ga^s\times\Gm^t$. Then $H$ is stationary if and only if $s\leq 1$.
\end{theorem}

Recall that we conjectured in \cite{YZ-Sta} that

\begin{Conj}\label{Conj-bad-nonsta}
Let $H$ be a $\dfg$ group definable in an $o$-minimal structure or a $p$-adically closed field. Then $H$ is bad if and only if $H$ is non-stationary.
\end{Conj}

Now we can prove Conjecture \ref{Conj-bad-nonsta} when $H$ is commutative.

\begin{coro}
Let $H$ be a commutative $\dfg$ group definable over $\Q$. Then $H$ is bad if and only if it is non-stationary.
\end{coro}
\begin{proof}
By Corollary \ref{H-product}, $H$ is virtually  $\Ga^s\times\Gm^t$ for some $s,t\in \N$. By Lemma \ref{equiv-stationarity-II}, $H$ is stationary iff $\Ga^s\times\Gm^t$ is stationary. If $H$ is bad, then by Remark \ref{bad-is-stationary} $H$ is non-stationary. conversely, if $H$ is non-stationary, then we see from Lemma \ref{GatimesGm^v-sta} that $s\geq 2$ and thus $\Ga^s\times\Gm^t$ is bad by Corollary \ref{u-geq-2}, so $H$ is bad.
\end{proof}



\section{The $\mu$-invariance}

Let $G$ be any group definable over $\Q$. We recall some notions from \cite{P-S: Top gup}.
The infinitesimal type of $G$ is the partial type over  $\Q$, denoted by $\mu_G(x)$ (or just by $\mu(x)$ if $G$ is clear), consisting of all formulas over $\Q$ defining an open neighborhood of $\text{id}_G$. To consistent with our earlier notation,  we identify $\mu (\M)$ (resp. $\mu_G(\M)$) with $\mu $  (resp. $\mu_G$).

\begin{notation}
    \begin{enumerate}
        \item If $\varphi(x)$ and $\psi(x)$ are $G$-formulas, then by $\varphi\cdot \psi$ we  denote the $G$-formula
        \[
        (\varphi\cdot \psi)(x)=\exists u\exists v (\varphi(u)\wedge\psi(v)\wedge x=uv)
        \]
        \item If $p(x)$ and $r(x)$ are (partial) $G$-types, then
        \[
        (p\cdot r)=\{(\varphi\cdot \psi)(x)|\ p\vdash \varphi(x),\ r\vdash \psi(x)\}
        \]
    \end{enumerate}
\end{notation}

\begin{definition}
Let $N$ be an elementary extension of $\Q$, $G$ a group definable over $\Q$, and $\mu$ the infinitesimal type of $G$.
\begin{enumerate}
    \item We say that $p\in S_G(N)$ is \emph{$\mu$-invariant} if $\mu \cdot p=p$;
    \item We say that $G$ is $\mu$-invariant if every global weakly generic type  on $G$ is $\mu$-invariant.
\end{enumerate}
\end{definition}

\begin{rmk}
It is easy to see that
\begin{enumerate}
    \item [(i)] $p\in S_G(N)$ is $\mu$-invariant if for any $a\models p$ and $\epsilon\models \mu $, we have that $ \epsilon a \models p$;
    \item  [(ii)] $G$ is $\mu$-invariant  iff every  weakly generic type on $G$ over every model $N\succ\Q$ is $\mu$-invariant invariant.
\end{enumerate}
\end{rmk}

We now fix $N\succ \Q$ as an arbitrary sufficiently saturated (small) submodel of $\M$.
\begin{definition}
 We say that $G$ has $\sharp$-property if for any $N^*\succ N$ ($|N^*|<|\M|$), each $H^{00}(N)$-invariant type $p\in S_G(N^*)$ is $\mu$-invariant.
\end{definition}

\begin{lemma}\label{lem-sharp}
Suppose that $G$ and $A$ are definably amenable groups definable over $\Q$.  Then we have:
\begin{enumerate}
    \item [(i)] If  $G$ admits $\sharp$-property, then it is $\mu$-invariant.
    \item [(ii)]  If $X$ is a finite index $\Q$-definable subgroup of $G$ and $X$ has $\sharp$-property, then $G$ also has $\sharp$-property.
    \item [(iii)] If $G^{00}=G^0$ and $G$ is virtually  $A$, then $G$  has $\sharp$-property iff $A$ has.
\end{enumerate}
\end{lemma}
\begin{proof}
    For (i), any global weakly generic type is $G^{00}$-invariant, thus is $\mu$-invariant by the definition of ``$\sharp$-property''.

    For (ii), take a $G^{00}(N)$-invariant type $p \in S_G(N^*)$, then $gp\in S_X(N^*)$ is also $G^{00}(N)$-invariant for any  $g\in G(\Q)$ since $gG^{00}=G^{00}g$. Take $g\in G(\Q)$ such that $gp\in S_X(N^*)$, then $gp$ is $\mu$-invariant. Since $g\mu=\mu g$, we see that $p$ is also $\mu$-invariant. This completes the proof.

    For (iii), let $X$ be a finite index $\Q$-definable subgroup of $G$ and  $f:X\rightarrow A$ a finite-to-one $\Q$-definable morphism such that $\img(f)$ has finite index in $A$. By part (ii), we may assume that $X=G$ and $A=\img(f)$.

     Let $N^*\succ N$ be a small submodel of $\M$, ${\cal G}\sq S_G(N^*)$ the space of $G^{0}(N)$-invariant types,  ${\cal A}\sq S_A(N^*)$ the space of $A^{0}(N)$-invariant types.  By Remark \ref{rmk-ker(f)}, we may assume that $\ker (f)\sq G^{0}(N)$, so $f$ is a one-one correspondence between the ${\cal G}$ and ${\cal A}$. Since $f(\mu_G)=\mu_A$, we see that $f(\mu_Gp)=\mu_Af(p)$ for any $p\in S_G(N^*)$. If $G$ has $\sharp$-property, then $\mu_Gp=p$ for any $p\in \cal G$, and thus $\mu_Aq=q$ for any $q\in\cal A$. So  $A$ has $\sharp$-property.

     Conversely, if $A$ has $\sharp$-property, then $\mu_Aq=q$ for any $q\in\cal A$.  By Lemma \ref{finite-to-one},  we have  $f^{-1}(\mu_A)=\mu_A\ker(f)$ and hence
    \[
   \mu_Gf^{-1}(q)=\mu_G\ker(f)f^{-1}(q)=f^{-1}(\mu_Aq)=f^{-1}(q)
    \]
for any $q\in\cal A$, so  $G$ has $\sharp$-property.
\end{proof}

\begin{lemma}\label{times-sharp}
Assume that $G_1$ and $ G_2$ are $\mathbb{Q}_p$ definable groups with $\sharp$-property. Then so is $G=G_1\times G_2$.
\end{lemma}
\begin{proof}
Let $N^*\succ N$ be a small submodel of $\M$ and $(\epsilon_1,\epsilon_2)\in\mu_G=\mu_{G_1}\times\mu_{G_2}$.
Let $h_1\in G_1$ and $h_2\in G_2$ such that $\tp(h_1,h_2/N^*)$ is $G^{00}(N)$-invariant.  Then we have that $\tp(h_1/ N^*, h_2)$ is $G_1^{00}(N)$-invariant. Since $G_1$ has $\sharp$-property, we have that
\[
\tp(\epsilon_1\cdot h_1/ N^*, h_2))=\tp(h_1/ N^*, h_2 ), \ \text{and hence}, \ \tp(\epsilon_1\cdot h_1,h_2/N^*)=\tp(h_1,h_2/N^*).
\]
Then $\tp(h_2/ N^*, \epsilon_1\cdot h_1 )$ is $G_2^{00}(N)$-invariant. A same argument shows  that
\[
\tp(\epsilon_1\cdot h_1,\epsilon\cdot h_2/N^*)=\tp(\epsilon_1\cdot h_1,h_2/N^*)
\]
and thus we have that
\[
\tp(\epsilon_1\cdot h_1,\epsilon\cdot h_2/N^*)=\tp(h_1, h_2/N^*)
\]
as required.
\end{proof}

Recall from \cite{Cluckers-Phd-Thesis} that
\begin{definition}
A  cell $A\sq \M$ is either a point or a set of the form
\[
\{t\in\M | \alpha \Box_1 v(t-c)\Box_2 \beta,\  (t-c)\in \lambda P_n(\M)\}
\]
with constants $n > 0$, $\lambda\in \mathbb Z$, $c\in \M$, $\alpha,\beta\in\Gamma_\M$, and $\Box_i$ either $<$ or no condition.
\end{definition}

\begin{fact}[\cite{Cluckers-Phd-Thesis}, Theorem 5.2.8.]\label{Cell-decom}
   Every definable set $X\sq \M$ is a finite disjoint union of the cells.
\end{fact}

\begin{lemma}\label{Ga_sharp}
$\Ga$ has $\sharp$-property.
\end{lemma}
\begin{proof}
Note that $\Ga^{00}(N)=\Ga(N)$. Let $\epsilon \in\mu_{\Ga}$, $N^*\succ N$ a small submodel of $  \M$, and $a\in\Ga $ with $\tp(a/N^*)$ a $\Ga(N)$-invariant type. By Fact \ref{Cell-decom}, to see that $\tp(\epsilon+a/N^*)=\tp(a/N^*)$, it suffices to show that for each $c\in\N^*$, $v(\epsilon+a-c)=v(a-c)$ and $P_n((\epsilon+a-c)/(a-c))$, for each $n\in \Z$.

Since $\tp(a/N^*)$ is $\Ga(N)$-invariant, we have that $v(a-c)<\Gamma_{N}$ for each $c\in N^*$. Hence, $v((\epsilon+a)-c)=v(\epsilon+(a-c))=v(a-c)$ for any $c\in N^*$. Also, for any  $c\in N^*$, we have that
\[
\frac{((\epsilon+a)-c)}{(a-c)}=1+\frac{\epsilon}{a-c}
\]
is infinitesimally close to $1$ over $N$, hence is an $n$-th power.

Hence, $\tp(\epsilon+a/N^*)=\tp(a/N^*)$ as required.
\end{proof}

\begin{lemma}\label{Gm_sharp}
$\Gm$ has $\sharp$-property.
\end{lemma}
\begin{proof}
Let $N^*\succ N$ and $a\in\Gm $ such that $\tp(a/N^*)$ is $\Gm^{00}(N)$-invariant. Firstly, $v(a)\geq v(a-c)$ for any $c\in N^*$.  Otherwise, there will be some $c\in N^*$ such that $v(a-c)>v(a)$, so $v(a)=v(c)$. Take $e\in \Gm^{00}(N)$ with $v(e)>0$. Then $v(ea)=v(c)$ as $\tp(a/N^*)$ is $\Gm^{00}(N)$-invariant. But $v(ea)=v(e)+v(a)=v(e)+v(c)>v(c)$, a contradiction.

We only need to show that $v(\delta a-c)=v(a-c)$ and $P_n((\delta a-c)/(a-c))$  for each $\delta \in\mu_{\Gm}$, $c\in N^*$, and  $n\in \Z$.
Take $\delta\in\mu_{\Gm}$, then $\delta=1+\epsilon$ for some $\epsilon\in \mu_{\Ga}$. For each $c\in N^*$, we have $v(\delta a-c)=v(\epsilon a+a-c)=v(a-c)$  since $v(\epsilon a)>v(a)\geq v(a-c)$.

Also, for any  $c\in \N^*$,
\[
\frac{\delta a-c}{a-c}=\frac{\epsilon a+a-c}{a-c}=1+\frac{\epsilon a}{a-c}.
\]
Since $v(\epsilon a/(a-c))\geq v(\epsilon a/a)=v(\epsilon)>\Z$, we have that $1+ (\epsilon a/(a-c))$ is an $n$-th power for each $n\in \N^{>0}$.
So $\tp(a/N^*)$ is $\mu_{\Gm}$-invariant as required.
\end{proof}

\begin{lemma}\label{mu-inv}
Let $H$ be a commutative dfg group over $\Q$, then $H$ has  $\sharp$-property.
\end{lemma}

\begin{proof}
By Corollary \ref{H-product}, $H$ is virtually $\Ga^s\times \Gm^t$ for some $s,t\in \N$.
We see from Lemma \ref{times-sharp}, \ref{Ga_sharp}, and \ref{Gm_sharp} that $\Ga^s\times \Gm^t$ has $\sharp$-property. By part (iii) of Lemma \ref{lem-sharp}, $H$  has $\sharp$-property.
\end{proof}

We see directly from Lemma \ref{mu-inv} that
\begin{coro}\label{mu-inv-II}
    Let $H$ be a commutative dfg group over $\Q$, then $H$ is  $\mu$-invariant.
\end{coro}

Let $M$ be a fixed model of $T=\pCF$ extending $\Q$ and $\M_0$ a monster model of $T^\ext_M$. It is clear that the reduct $\M_0\upharpoonright_L$ of  $\M_0$ to $L$ is also saturated and homogeneous, for convenience, we assume that  $\M_0\upharpoonright_L=\M$.
It is easy to see that the partial type $\mu_H$ is the same whether computed in $\M$ of $T$ or $\M_0$ of $T^\ext_M$.
\begin{lemma}\label{mu-inv-ext}
Let $H$ be a commutative $\dfg$ group definable over $\Q$. Then  $H$ is also  $\mu$-invariant when we compute it in $T^\ext_M$.
\end{lemma}
\begin{proof}
Note that $T^\ext_M$ also has $\NIP$. By Fact \ref{G00=G00ext}, $H^{00}$ is the same whether computed in $T$ or $T^\ext_M$. Let $N_0\succ M^\ext$ be a sufficiently saturated model and $N$  a reduct of $N_0$ to $L$. Let $p\in S_H(N_0)$ be a weakly generic type. By Fact \ref{G00-inv=WG}, $p$ is $H^{00}(N)$-invariant.    By Fact \ref{Sim-Exercise}, there is $N^*\prec \M$  such that $N\prec N^*$ and  $S_H(N_0)$ is naturally homeomorphic to $S_{H, N}(N^*)$.  We consider  $p$ as an element of $S_{H, N}(N^*)\subset S_H(N^*)$, then $p$ is $H^{00}(N)$-invariant and we see from
 Lemma \ref{mu-inv}  that $p$ is   $\mu$-invariant.
\end{proof}

\section{The main theorem}

In this section, we fix $M$ as an elementary extension of $\Q$.
We denote $M^\ext$ by $M_0$. By $\M_0$,  Let $\M_0\succ M_0$ be a monster model of $T^\ext_M$. Also, we assume that $\M$, the reduct of $\M_0$ in $L$, is a monster model of $T=\pCF$ extending $M$.

Also, we fix $G$ as a commutative group definable over $\Q$ which is not a $\dfg$.
Note that by Fact \ref{dfg-fsg-decomp}, $G$ admits a $\Q$-definable short exact sequence
\[
1\rightarrow H\rightarrow G\rightarrow_\pi C\rightarrow 1,
\]
where $C$ is a $\fsg$ group with $\dim(C)\geq 1$ and $H$ a $\dfg$ group. By Fact \ref{fact-fsg-ext} and Fact \ref{dfg-fact-2}, $H$ and $C$ also have $\dfg$ and $\fsg$, respectively, when we compute them in $T_M^\ext$.

 Let $f:C\rightarrow G$ be a $\Q$-definable section of $\pi$. Then any $g\in G$ can be written uniquely as $f(c)h$  for $c=\pi(g)\in C$ and $h=(f(c))^{-1}g\in H$. Moreover, we can assume that $\img(f)$ is contained in a definably compact subset of $G$: Let $U(\Q)$ be an open compact subgroup of $G(\Q)$ definable over $\Q$.  Then $\pi(U(\Q))$ is an open subgroup of $C(\Q)$. Since $C(\Q)$ is compact, we can find  $c_1,...,c_n\in C(\Q)$ with $n\in \N$, such that $\bigcup_{i\leq n}c_i\cdot\pi(U(\Q))=C(\Q)$. Let $g_i\in G(\Q)$ such that $\pi(g_i)=c_i$ for $i=1,...,n$ and $V(\Q)= U(\Q)\cup \bigcup_{i\leq n} g_i\cdot U(\Q)$. We have that $V(\Q)$ is  definable over $\Q$ and $\pi(V(\Q))=C(\Q)$. Now $V(\Q)$ is an open compact neighbourhood of $1_G$ as it is a finite union of open compact sets.  Let $V=V(\M)$, by definable Skolem functions, we can find a $\Q$-definable section $f$ of $\pi\upharpoonright_{V}$. Hence, we can assume that $f(C)\subset V$. Replacing $V$ by $V\cup V^{-1}$, we may assume that $V=V^{-1}$.
Let $\eta: C\times C\rightarrow H$ defined by $\eta(c_1,c_2)=f(c_1c_2)^{-1}f(c_1)f(c_2)$.

\begin{lemma}\label{eta-mu}
For any $c_1,c_2\in C$, there is $h_0\in H(\Q)$ such that $\eta(c_1,c_2)\in\mu_H h_0$.
\end{lemma}
\begin{proof}
Firstly, we show that $H$ is closed in $G$. It is clear that $\cl(H)$ is a closed subgroup of $G$ definable over $\Q$. By Denef's cell decomposition in \cite{Denef-II}, we see that $\dim(H)=\dim(\cl(H))$. Then $H$ is a subgroup of $\cl(H)$ with nonempty interior, which implies that $H$ is an open subgroup of $\cl(H)$. Thus, $H$ is a closed subgroup of $\cl(H)$. Hence, $H=\cl(H)$ is closed in $G$.

Let $W=V\cdot V\cdot V$, then $W$ is definably compact as $V$ is. It is clear that $\eta(C)$ is contained in $W\cap H$. Since $W(\Q)$ is compact and $H(\Q)$ is closed in $G(\Q)$, we have that $W(\Q)\cap H(\Q)$ is compact. Then for any $c_1,c_2\in C$, $\eta(c_1,c_2)\in (W\cap H)$ has a standard part $h_0\in (W\cap H)(\Q)$, namely, $\eta(c_1,c_2)\in(\mu_G h_0)\cap H=\mu_H h_0$.
\end{proof}

\subsection{Local case}
We now describe the almost periodic types in $S_G(M_0)$.
Let $\I$ be the space of generic types in $S_C(M_0)$, then $\I=\AP(S_C(M_0))$ is the unique minimal subflow of $S_C(M_0)$ and it is also a bi-ideal of the semigroup $(S_C(M_0),*)$ (see Fact \ref{AP-bi-ideal}). By $f(\I)$ we denote the set $\{f(q)|\ q\in \I\}$. Let $\J=\AP(S_H(M_0))$ be the collection of all almost periodic types   in $S_H(M_0)$ which is also the union of all minimal subflows of $S_H(M_0)$.

\begin{fact}[\cite{YZ-Sta}]\label{Min-flow=I-J}
$f(\I)*\J\subseteq \AP(S_G(M_0))$.
\end{fact}

\begin{fact}[\cite{YZ-Sta}]\label{AP=q-p-r}
Let $r\in S_G(M_0)$. Then $r$ is almost periodic iff $r=f(q)*p*r$ for some  $q\in \I$ and  $p\in \J$.
\end{fact}

Note that by Theorem 5.1 of  \cite{Boxall-Kestner} that $T_M^\ext=\Th(M_0)$ is also a distal theory. By Lemma 2.16 in \cite{Simon-distal}, we have

\begin{fact}\label{oth-ext}
    For any $N\models T^\ext$, if and $N'\succ N$ is $|N|^+$-saturated, $p(x)\in S(N')$ definable over $N$, and $q(y)\in S(N')$ finitely satisfiable in $N$. Then 
    $p(x)\cup q(y)$ implies a complete  $(x,y)$-type over $N'$. In fact, if $a\models p$ and $b\models q$, then $\tp(a/N',b)$ is the unique heir of $\tp(a/N)$ and $\tp(b/N',a)$ is finitely satisfiable in $N$.
\end{fact}

\begin{fact}\label{dfg-fact-2}\cite{CPS-Ext}
Let $N$ be a sufficiently saturated extension of $M^\ext$. Then $p\in S_H(M_0)$ is almost periodic iff its unique  heir $\bar p$ over $N$ is weakly generic.
\end{fact}

\begin{lemma}\label{AP=f(I)J}
$\AP(S_G(M_0))=f(\I)\ast \J$.
\end{lemma}

\begin{proof}
By fact \ref{Min-flow=I-J}, we only need to prove that $\AP(S_G(M_0))\subseteq f(\I)\ast \J$.

Let $r\in S_G(M_0)$ be almost periodic, then by Fact \ref{AP=q-p-r} we have $r=f(q)*p*r$ for some $q\in \I$ and $p\in \J$. Let $N$ be a sufficiently  extension of $M_0$. Let $h\in H(N)$ realize the unique heir of $p$ over $N$,   Let $g^*=f(c^*)h^*\in G$ realize the unique heir of $r$ over $\dcl(N,h)$ with $c^*=\pi(g^*)\in C$ and $h^*\in H$. By Fact \ref{fact-fsg}, $q$ has a unique generic  extension $\bar q$ over $\dcl(N,h,g^*)$, which is finitely satisfiable in $M_0$. Let $c\in C$ realize  $\bar q$,  then
\[
f(q)*p*r=\tp(f(c)hf(c^*)h^*/M_0)=\tp(f(c)f(c^*)hh^*/M_0)=\tp(f(cc^*)\eta(c,c^*)hh^*/M_0).
\]
Clearly, $\tp(h^*/M_0,h)$ is the unique heir of $\tp(h^*/M_0)$ since $h^*\in\dcl(M_0,g^*)$. We have that $\tp(hh^*/M_0)=p\ast \tp(h^*/M_0)$. By Fact \ref{AP-bi-ideal}, we have that $\tp(hh^*/M_0)\in\J$ since $p\in\J$. 
Since both $\tp(h/N)$  and $\tp(h^*/N,h)$ are definable over $M_0$, we have  that $\tp(hh^*/N)$ is definable over $M_0$, and thus the heir extension of $\tp(hh^*/M_0)$. We see from  Fact \ref{dfg-fact-2} that $\tp(hh^*/N)$ is a weakly generic type.

By Lemma \ref{eta-mu} , there is $h_0\in H(\Q)$ such that $\eta(c,c^*)=\epsilon h_0$ for some $\epsilon\in \mu_H$, and hence, by Lemma \ref{mu-inv-ext},
\[
\tp(\eta(c,c^*)hh^*/N)=\tp(\epsilon h_0hh^*/N)=h_0\cdot\tp(hh^*/N).
\]
So $\tp(\eta(c,c^*)hh^*/N)$ is also definable over $M_0$.
As $\tp(cc^*/N)=\bar q* \tp(c^*/N)$ is a generic type on $C$, thus is a  coheir extension of $\tp(cc^*/M_0 \epsilon h_0hh^*)$. By Fact \ref{oth-ext}, $\tp(f(cc^*)/M_0,\eta(c,c^*)hh^*)$ is finitely satisfiable in $M_0$, so we have
\[
 r=   f(q)*p*r=\tp(f(cc^*)\eta(c,c^*)hh^*/M_0)=\tp(f(cc^*)/M_0)\ast\tp(\eta(c,c^*)hh^*/M_0).
\]
Since $\tp(cc^*/M_0)\in \I$ and $\tp(\eta(c,c^*)hh^*/M_0)=h_0\tp(hh^*/M_0)\in\J$, we have that $r\in f(\I)\ast\J$ as required.
\end{proof}

\begin{fact}\cite{YZ-Sta}\label{heir-neq-coheir}
Let $p=\tp(e/\Q)\in S_1(\Q)$ be a non-algebraic type. Suppose that $p_1$ is the unique heir  $p$ over $\Q,e$, then $p_1$ is not a finitely satisfiable in $\Q$.
\end{fact}

\begin{fact}\cite{YZ-Sta}\label{construct-strongly-f-generic}
 Let $N$ be an $|M_0|^+$-saturated extension of $M_0$, $\tp(c^*/N )$   a generic type on $C$, and $\tp(h^*/ N, c^* )$  a strongly $f$-generic type on $H$ over $\Q$. Then $\tp(f(c^*)h^*/N)$ is a strongly $f$-generic type on $G$ over $\Q$.
\end{fact}

\begin{lemma}\label{Local-bad-neg}
Let $H$ be a bad $\dfg$ group (in $T$). Then $\WGen(S_G(M_0))$ is a proper subset of $\AP(S_G(M_0))$.
\end{lemma}
\begin{proof}
Let $N_0$ be a sufficiently saturated extension of $M_0$, and $N=N_0\upharpoonright_L$ the reduct of $N_0$ to $L$. Let $\tp(c^*/N)$ be a generic type on $C$. By badness of $H$, we can find $\tp(h^*/ N,c^*)$, a strongly $f$-generic type on $H$ over $\Q$, and a $\Q$-definable function $\theta$, such that $\tp(\theta(h^*)/N,c^*)$ is finitely satisfiable in $\Q$. We see from Fact \ref{construct-strongly-f-generic} that $p=\tp(f(c^*)h^*/N)\in S_G(N)$ is strongly $f$-generic over $\Q$, hence is $G^{00}(N)$-invariant.

Let us consider $p$ as a partial type over $N_0$ (in $T_M^\ext$). By Fact \ref{G00=G00ext}, $G^{00}$ is the same whether computed it in $T$ or $T^{\ext}_M$, so   $p$ is also  weakly generic  in $T_M^\ext$. Let $p^*\in S_G(N_0)$ be any weakly generic extension of $p$ and $p_0\in S_G(M_0)$  the restriction  of $p^*$ to $M_0$, then $p_0$ is weakly generic.

Suppose for a contradiction that $p_0$ is almost periodic. Let $q=p_0\upharpoonright_{L}\in S_G(M)$ be the restriction of $p_0$ to the language $L$, and $q_0=q\upharpoonright_{\Q}\in S_G(\Q)$   the restriction of $q$  to $\Q$  (in $T$). Then $q$ is almost periodic in $S_G(M)$ by part (i) of Fact \ref{min-flow-restriction}, and thus $q_0$ is almost periodic in $S_G(\Q)$ by part (ii) of Fact \ref{min-flow-restriction}.
Clearly, $q_0=\tp(f(c^{*})h^{*}/\Q)$. Applying Lemma \ref{AP=f(I)J} to $\Q^\ext=\Q$, we have that $\tp(f(c^*)/\Q,h^*)$ is finitely satisfiable in $\Q$. We conclude that both $\tp(f(c^*)/\Q,\theta(h^*))$ and $\tp(\theta(h^*)/\Q,f(c^*))$ are finitely satisfiable in $\Q$, which contradicts to Fact \ref{heir-neq-coheir}.
\end{proof}

Then, combining Lemma \ref{Local-bad-neg} with Fact \ref{sta-AP=WG},
we have that
\begin{theorem}
$\AP(S_G(M_0))=\WGen(S_G(M_0))$ if and only if $G$ is stationary.
\end{theorem}

\subsection{Global case}

Now we consider the global case.

\begin{lemma}\label{Global-bad-neg}
Let $H$ be a bad $\dfg$ group. Then $S_G(\M)$ has a weakly generic type which is not almost periodic.
\end{lemma}
\begin{proof}
Let $c^*$ realize a generic type  $p\in S_C(\M)$. By badness of $H$, we can find a strongly $f$-generic type $\tp(h^*/ \M, c^*)$ on $H$   and a $\Q$-definable function $\theta$, such that $\theta(h^*)$ is finitely satisfiable in $\Q$. We see from Fact \ref{construct-strongly-f-generic} that $\tp(f(c^*)h^*/\M)$ is strongly $f$-generic on $G$ over $\Q$.

Suppose for a contradiction that $\tp(f(c^*)h^*/\M)$ is almost periodic.  Then $\tp(f(c^*)h^*/\Q)$ is almost periodic in $S_G(\Q)$ by part (ii) of Fact \ref{min-flow-restriction}. By Lemma \ref{AP=f(I)J}, we have that $\tp(c^*/\Q,h^*)$ is finitely satisfiable in $\Q$. We conclude that both $\tp(c^*/\Q,\theta(h^*))$ and $\tp(\theta(h^*)/\Q,c^*)$ are finitely satisfiable in $\Q$, which contradicts to Fact \ref{heir-neq-coheir}.
\end{proof}

Then, combining Lemma \ref{Global-bad-neg} with Fact \ref{sta-AP=WG},
we have that
\begin{theorem}
$\AP(S_G(\M))=\WGen(S_G(\M))$ if and only if $G$ is stationary.
\end{theorem}

\end{document}